\documentclass[11pt]{amsart}
\usepackage{amsmath,amsthm,amssymb,amscd}

\theoremstyle{plain}

\newtheorem{thm}{Theorem}[section]
\newtheorem{proposition}[thm]{Proposition}

\newtheorem{lem}[thm]{Lemma}

\newtheorem{prop}[thm]{Proposition}

\theoremstyle{definition}

\newtheorem{defn}[thm]{Definition}
\newtheorem{ex}[thm]{Example}
\newtheorem{rem}[thm]{Remark}

\newcommand\cA{{\mathcal{A}}}
\newcommand\cB{{\mathcal{B}}}

\newcommand\cD{{\mathcal{D}}}

\newcommand\cF{{\mathcal{F}}}
\newcommand\cG{{\mathcal{G}}}
\newcommand\cH{{\mathcal{H}}}

\newcommand\cK{{\mathcal{K}}}
\newcommand\cL{{\mathcal{L}}}

\newcommand\RR{\mathbb{R}}

\newcommand\CC{\mathbb{C}}

\newcommand\NN{\mathbb{N}}

\newcommand\Ker{\mathop{\rm Ker}\nolimits}
\newcommand\Dom{\mathop{\rm Dom}\nolimits}
\newcommand\supp{\mathop{\rm supp}\nolimits}
\newcommand{\cinf}{C^\infty}

\begin{document}
\title[Egorov's theorem for transversally elliptic operators]{Egorov's theorem for
transversally elliptic operators on foliated manifolds and
noncommutative geodesic flow}
\author{Yuri A. Kordyukov}
\address{Institute of Mathematics, Russian Academy of Sciences,
Ufa, Russia} \email{yuri@imat.rb.ru}

\begin{abstract}
The main result of the paper is Egorov's theorem for
trans\-ver\-sal\-ly elliptic operators on compact foliated
manifolds. This theorem is applied to describe the noncommutative
geodesic flow in noncommutative geometry of Riemannian foliations.
\end{abstract}
\maketitle
\hyphenation{trans-ver-sal-ly}

\bibliographystyle{plain}

\section*{Introduction}
Egorov's theorem \cite{egorov} is one of the fundamental results
in microlocal analysis that relates the quantum evolution of
pseudodifferential operators with the classical dynamics of
principal symbols.

Let $P$ be a positive, self-adjoint, elliptic, first order
pseudodifferential operator on a compact manifold $M$ with the
positive principal symbol $p\in S^1(T^*M\setminus 0)$. Let $f_t$
be the bicharacteristic flow of the operator $P$, that is, the
Hamiltonian flow of $p$ on $T^*M$. For instance, one can consider
$P=\sqrt{\Delta_M}$, where $\Delta_M$ is the Laplace operator of a
Riemannian metric $g_M$ on $M$. Then the bicharacteristic flow of
the operator $P$ is the geodesic flow of the metric $g_M$.

Egorov's theorem states that, for any pseudodifferential operator
$A$ of order $0$ with the principal symbol $a\in S^0(T^*M\setminus
0)$, the operator $A(t)=e^{itP}Ae^{-itP}$ is a pseudodifferential
operator of order $0$. The principal symbol $a_t\in
S^0(T^*M\setminus 0)$ of this operator is given by the formula $$
a_t(x,\xi)=a(f_t(x,\xi)), \quad (x,\xi)\in T^*M\setminus 0. $$

The main result of this paper is a version of Egorov's theorem for
transversally elliptic operators on compact foliated manifolds.
This theorem is applied to describe the noncommutative geodesic
flow in noncommutative geometry of Riemannian foliations.

\newpage
\centerline{\textbf{Table of Contents}}
\bigskip

\noindent 1. Preliminaries and main results

1.1 Transverse pseudodifferential calculus

1.2 Transverse bicharacteristic flow

1.3 Egorov's theorem

1.4 Noncommutative geodesic flow on foliated manifolds

\noindent 2. Proof of the main theorem

2.1 The case of elliptic operator

2.2 The general case

\noindent 3 Noncommutative geometry of foliations

\section{Preliminaries and main results}
\label{Eg:sect}
\subsection{Transverse pseudodifferential calculus}\label{s:trpdo}
Throughout in the paper, $(M,{\mathcal F})$ is a compact foliated
manifold, $E$ is a Hermitian vector bundle on $M$,
$\operatorname{dim} M=n, \operatorname{dim} \cF=p, p+q=n$.

We will consider pseudodifferential operators, acting on
half-den\-si\-ti\-es. For any vector bundle $V$ on $M$, denote by
$|V|^{1/2}$ the associated half-density vector bundle. Let
$C^{\infty}(M,E)$ denote the space of smooth sections of the
vector bundle $E\otimes |TM|^{1/2}$, $L^2(M,E)$ the Hilbert space
of square integrable sections of $E\otimes |TM|^{1/2}$,
${\cD}'(M,E)$ the space of distributional sections of $E\otimes
|TM|^{1/2}$, ${\cD}'(M,E)=C^{\infty}(M,E)'$, and $H^s(M,E)$ the
Sobolev space of order $s$ of sections of $E\otimes |TM|^{1/2}$.
Finally, let $\Psi^{m}(M,E)$ denote the standard classes of
pseudodifferential operators, acting in $C^{\infty}(M,E)$.

We will use the classes $\Psi^{m,-\infty}(M,{\mathcal F},E)$ of
transversal pseudodifferential operators. Let us briefly recall
its definition, referring the reader to \cite{noncom} for more
details.

We will consider foliated coordinate charts $\varkappa : U\subset
M\stackrel{\sim}{\longrightarrow} I^{n}$ on $M$ with coordinates
$(x,y)\in I^p\times I^q$ ($I$ is the open interval $(0,1)$) such
that the restriction of $\cF$ to $U$ is given by the sets $y={\rm
const}$. We will always assume that foliated charts are regular.
Recall that a foliated coordinate chart $\varkappa : U\subset
M\stackrel{\sim}{\longrightarrow} I^{n}$ is called regular, if it
admits an extension to a foliated coordinate chart
$\bar{\varkappa} : V\subset M
\stackrel{\sim}{\longrightarrow}(-2,2)^{n}$ with $\bar{U}\subset
V$.

A map $f:U\subset M\rightarrow \RR^q$ is called a distinguished
map, if $f$ locally has the form $pr_{nq}\circ \varkappa$, where
$\varkappa : V\subset U\stackrel{\sim}{\longrightarrow} I^{n}$ is
a foliated chart and $pr_{nq}:\RR^n=\RR^p\times\RR^q \rightarrow
\RR^q$ is the natural projection. Let $D_x$ denote the set of
germs of distinguished maps from $M$ to $\RR^q$ at a point $x\in
M$. For any leafwise continuous curve $\gamma$ from $x$ to $y$,
let $h_\gamma: D_x\rightarrow D_y$ be the holonomy map associated
with $\gamma$. This is the generalization of Poincar\'e's first
return map from flows to foliations.

Let $\varkappa: U\to I^p\times I^q, \varkappa': U'\to I^p\times
I^q$, be two foliated charts, $\pi=pr_{nq}\circ \varkappa: U\to
\RR^q$, $\pi'=pr_{nq}\circ \varkappa': U'\to \RR^q$ the
corresponding distinguished maps. The foliation charts
$\varkappa$, $\varkappa'$ are called compatible, if, for any $m\in
U$ and $m'\in U'$ such that $m=\varkappa^{-1}(x,y)$,
$m'={\varkappa'}^{-1}(x',y)$ with the same $y$, there is a
leafwise path $\gamma$ from $m$ to $m'$ such that the
corresponding holonomy map $h_{\gamma}$ takes the germ $\pi_m$ of
the map $\pi$ at $m$ to the germ $\pi'_{m'}$ of the map $\pi'$ at
$m'$.

Let $\varkappa: U\subset M\rightarrow I^p\times I^q, \varkappa':
U' \subset M\rightarrow I^p\times I^q$, be two compatible foliated
charts on $M$ equipped with trivializations of the vector bundle
$E$ over them. Consider an operator $A:C^{\infty}_c(U,\left.
E\right|_U)\to C^{\infty}_c(U',\left. E\right|_{U'})$ given in the
local coordinates by the formula
\begin{equation}\label{loc}
Au(x,y)=(2\pi)^{-q} \int e^{i(y-y')\eta}k(x,x',y,\eta)
u(x',y') \,dx'\,dy'\,d\eta,
\end{equation}
where $k \in S ^{m} (I ^{p} \times I^p\times I^q\times {\RR}^{q},
{\cL}({\CC}^r))$, $u \in C^{\infty}_{c}(I^{n}, {\CC}^r), x \in
I^{p}, y \in I^{q}$ with the Schwartz kernel, compactly supported
in $U\times U'$ (here $r=\operatorname{rank} E$).

Recall that a function $k\in C^{\infty}(I^{p} \times I^{p} \times
I^{q} \times {\RR}^{q},{\mathcal L}({\CC}^r))$ belongs  to  the
class $S^{m}(I^{p}\times I^{p}\times I^{q} \times {\RR}^{q},
{\cL}({\CC}^r))$, if, for any multiindices $\alpha $ and $\beta $,
there exists a constant $C_{\alpha \beta} > 0$ such that
\begin{multline*}
|\partial^{\alpha}_{\eta} \partial^{\beta}_{(x,x',y)}k(x,x',y,\eta
)| \leq C_{\alpha \beta }(1 +\vert \eta \vert )^{ m-\vert \alpha
\vert },\\ (x,x',y)\in I^{p}\times I^p\times I^q,\quad \eta \in
{\RR}^{q}.
\end{multline*}
We will consider only classical symbols $k$, which
can be represented as an asymptotic sum
$
k(x,x',y,\eta)\sim \sum_{j=0}^{\infty} \theta(\eta)
k_{z-j}(x,x',y,\eta),
$
where $k_{z-j}\in C^{\infty}(I^{p}\times I^{p}\times I^{q} \times
({\RR}^{q}\backslash \{0\}), {\cL}({\CC}^r))$ is homogeneous in
$\eta$ of degree $z-j$, and $\theta$ is a smooth function on
${\RR}^{q}$ such that $\theta(\eta)=0$ for $|\eta|\leq 1$,
$\theta(\eta)=1$ for $|\eta|\geq 2$.

The operator $A$ extends to an operator in $C^\infty(M,E)$ in a
trivial way. The resulting operator is called an elementary
operator of class $\Psi ^{m,-\infty}(M,{\mathcal F},E)$.

The  class $\Psi ^{m,-\infty}(M,{\mathcal F},E)$ consists of all
operators $A$ in $C^{\infty}(M,E)$, which can be represented in
the form $A=\sum_{i=1}^k A_i + K$, where $A_i$ are elementary
operators of class $\Psi ^{m,-\infty}(M,{\mathcal F},E)$,
corresponding to some pairs $\varkappa_i,\varkappa'_i$ of
compatible foliated charts, $K\in \Psi ^{-\infty}(M,E)$. Put
$\Psi^{*,-\infty}(M,{\mathcal
F},E)=\bigcup_m\Psi^{m,-\infty}(M,{\mathcal F},E)$.

Let $G$ be the holonomy groupoid of ${\cF}$. We will  briefly
recall its definition. Let $\sim_h$ be the equivalence relation on
the set of continuous leafwise paths $\gamma:[0,1] \rightarrow M$,
setting $\gamma_1\sim_h \gamma_2$ if $\gamma_1$ and $\gamma_2$
have the same initial and final points and the same holonomy maps.
The holonomy groupoid $G$ is the set of $\sim_h$ equivalence
classes of continuous leafwise paths. $G$ is equipped with the
source and the range maps $s,r:G\rightarrow M$ defined by
$s(\gamma)=\gamma(0)$ and $r(\gamma)=\gamma(1)$. We will identify
a point $x\in  M$ with the element of $G$ given by the
corresponding constant path: $\gamma(t)=x, t\in [0,1]$. Recall
also that, for any $x\in M$, the set $G^x=\{\gamma\in G :
r(\gamma)=x\}$ is the covering of the leaf through the point $x$
associated with the holonomy group $G^x_x$ of this leaf,
$G^x_x=\{\gamma\in G : s(\gamma)=x, r(\gamma)=x\}$.

Any pair of compatible foliated charts $\varkappa: U\to I^p\times
I^q, \varkappa': U'\to I^p\times I^q$ defines a foliated chart
$V\to I^p\times I^p\times I^q$ on $G$ as follows. The coordinate
patch $V$ consists of all $\gamma\in G$ from
$m=\varkappa^{-1}(x,y)\in U$ to $m'={\varkappa'}^{-1}(x',y)\in U'$
such that the corresponding holonomy map $h_{\gamma}$ takes the
germ $\pi_m$ of the distinguished map $\pi=pr_{nq}\circ\varkappa$
at $m$ to the germ $\pi'_{m'}$ of the distinguished map
$\pi'=pr_{nq}\circ\varkappa'$ at $m'$, and the coordinate map
takes such a $\gamma$ to $(x,x',y)\in I^p\times I^p\times I^q$.

Denote by $N^*{\mathcal F}$ the conormal bundle to ${\mathcal F}$.
For any $\gamma\in G, s(\gamma)=x, r(\gamma)=y$, the
codifferential of the corresponding holonomy map defines a linear
map $dh^*_{\gamma}: N^*_y{\mathcal F}\to N^*_x{\mathcal F}$. Let
${\mathcal F}_N$ be the linearized foliation in
$\tilde{N}^*{\mathcal F}=N^*{\mathcal F}\setminus 0$ (cf., for
instance, \cite{Molino}). The leaf of the foliation $\cF_N$
through $\nu\in \tilde{N}^*{\mathcal F}$ is the set of all points
$dh_{\gamma}^{*}(\nu)\in \tilde{N}^*{\mathcal F}$, where
$\gamma\in G, r(\gamma)=\pi(\nu)$ (here $\pi :T^*M\to M$ is the
bundle map). The leaves of the foliation $\cF_N$ have trivial
holonomy. Therefore, the holonomy groupoid $G_{{\mathcal F}_N}$ of
${\mathcal F}_N$ consists of all pairs $(\gamma,\nu)\in G\times
\tilde{N}^*{\mathcal F}$ such that $r(\gamma)=\pi(\nu)$ with the
source map $s_N:G_{{\mathcal F}_N}\rightarrow \tilde{N}^*{\mathcal
F}, s_N(\gamma,\nu)=dh_{\gamma}^{*}(\nu)$ and the range map
$r_N:G_{{\mathcal F}_N}\rightarrow \tilde{N}^*{\mathcal F},
r_N(\gamma,\nu)=\nu$. We have a map $\pi_G:G_{{\mathcal
F}_N}\rightarrow G$ given by $\pi_G(\gamma,\nu)=\gamma$.

Denote by $\pi^*E$ the lift of the vector bundle $E$ to
$\tilde{N}^*\cF$ via the bundle map $\pi:\tilde{N}^*\cF\to M$ and
by ${\mathcal L}(\pi^*E)$ the vector bundle on $G_{\cF_N}$, whose
fiber at a point $(\gamma,\nu)\in G_{\cF_N}$ is the space
$\cL((\pi^*E)_{s_N(\gamma,\nu)}, (\pi^*E)_{r_N(\gamma,\nu)})$ of
linear maps from $(\pi^*E)_{s_N(\gamma,\nu)}$ to
$(\pi^*E)_{r_N(\gamma,\nu)}$. There is a natural foliation $\cG_N$
on $G_{{\mathcal F}_N}$. The leaf of $\cG_N$ through a point
$(\gamma,\nu)\in G_{{\mathcal F}_N}$ is the set of all $(\gamma',
\nu')\in G_{{\mathcal F}_N}$ such that $\nu$ and $\nu'$ lie in the
same leaf in $\cF_N$. Let $|T{\mathcal G}_N|^{1/2}$ be the line
bundle of leafwise half-densities on $G_{{\mathcal F}_N}$ with
respect to the foliation ${\mathcal G}_N$. It is easy to see that
$$ |T{\mathcal G}_N|^{1/2}=r_N^*(|T{\mathcal F}_N|^{1/2})\otimes
s_N^*(|T{\mathcal F}_N|^{1/2}), $$ where $s_N^*(|T{\mathcal
F}_N|^{1/2})$ and $r_N^*(|T{\mathcal F}_N|^{1/2})$ denote the
lifts of the line bundle $|T{\mathcal F}_N|^{1/2}$ of leafwise
half-densities on $N^*{\mathcal F}$ via the source and the range
mappings $s_N$ and $r_N$ respectively.

A section $k\in C^{\infty}(G_{{\mathcal F}_N}, {\mathcal
L}(\pi^*E)\otimes |T{\mathcal G}_N|^{1/2})$ is said to be properly
supported, if the restriction of the map $r:G_{\cF_N}\to
\tilde{N}^*\cF$ to $\supp k$ is a proper map. Consider the space
$C^{\infty}_{prop}(G_{{\mathcal F}_N}, {\mathcal L}(\pi^*E)\otimes
|T{\mathcal G}_N|^{1/2})$ of smooth, properly supported sections
of ${\mathcal L}(\pi^*E)\otimes |T{\mathcal G}_N|^{1/2}$. One can
introduce the structure of involutive algebra on
$C^{\infty}_{prop}(G_{{\mathcal F}_N}, {\mathcal L}(\pi^*E)\otimes
|T{\mathcal G}_N|^{1/2})$ by the standard formulas (cf.
(\ref{e:a})). Let $S^{m}(G_{{\mathcal F}_N},{\mathcal
L}(\pi^*E)\otimes |T{\mathcal G}_N|^{1/2})$ be the space of all
$s\in C^{\infty}_{prop}(G_{{\mathcal F}_N},{\mathcal
L}(\pi^*E)\otimes |T{\mathcal G}_N|^{1/2})$ homogeneous of degree
$m$ with respect to the action of $\RR$ given by the
multiplication in the fibers of the vector bundle
$\pi_G:G_{{\mathcal F}_N}\rightarrow G $. By \cite{noncom}, there
is the half-density principal symbol mapping
\begin{equation}\label{e:symbol}
\sigma: \Psi^{m,-\infty}(M,{\mathcal F},E)\rightarrow
S^m(G_{{\mathcal F}_N},{\mathcal L}(\pi^*E)\otimes |T{\mathcal
G}_N|^{1/2}),
\end{equation}
which satisfies
\[
\sigma_{m_1+m_2}(AB)=\sigma_{m_1}(A)\sigma_{m_2}(B),\quad
\sigma_{m_1}(A^*)=\sigma_{m_1}(A)^*
\]
for any $A\in \Psi^{m_1,-\infty}(M,{\mathcal F},E)$ and $B\in
\Psi^{m_2,-\infty}(M,{\mathcal F},E)$.

\begin{ex}
Consider a foliated coordinate chart $\varkappa : U\subset
M\stackrel{\sim}{\longrightarrow} I^{n}$ on $M$ with coordinates
$(x,y)\in I^p\times I^q$. One has the corresponding coordinate
chart in $T^*M$ with coordinates given by $(x,y,\xi,\eta)\in
I^p\times I^q\times \RR^p\times \RR^q$. In these coordinates, the
restriction of the conormal bundle $N^*{\mathcal  F}$ to $U$ is
given by the equation $\xi=0$. So we have a coordinate chart
$\varkappa_n : U_1\subset N^*{\mathcal  F}
\stackrel{\sim}{\longrightarrow} I^p\times I^q\times \RR^q$ on
$N^*{\mathcal  F}$ with the coordinates $(x,y,\eta)\in I^p\times
I^q\times \RR^q$. The coordinate chart $\varkappa_n$ is a foliated
coordinate chart for the linearized foliation ${\mathcal  F}_N$,
and the restriction of ${\mathcal  F}_N$ to $U_1$ is given by the
level sets $y= {\rm const}, \eta={\rm const}$.

Now let $\varkappa: U\subset M\rightarrow I^p\times I^q,
\varkappa': U' \subset M\rightarrow I^p\times I^q$, be two
compatible foliated charts on $M$. Then the corresponding foliated
charts $\varkappa_n: U_1\subset N^*{\mathcal  F}\rightarrow
I^p\times I^q \times \RR^q, \varkappa'_n: U'_1 \subset
N^*{\mathcal  F} \rightarrow I^p\times I^q \times \RR^q,$ are
compatible with respect to the foliation ${\mathcal  F}_N$. So
they define a foliated chart $V$ on the foliated manifold
$(G_{{\mathcal F}_N},{\mathcal G}_N)$ with the coordinates
$(x,x',y,\eta) \in I^p\times I^p\times I^q \times \RR^q$, and the
restriction of ${\mathcal  G}_N$ to $V$ is given by the level sets
$y= {\rm const}, \eta={\rm const}$. The principal symbol
$\sigma_m(A)$ of an operator $A$ given by the formula (\ref{loc})
is the half-density $k_m(x,x',y,\eta)\, |dx|^{1/2}\,|dx'|^{1/2}$,
where $k_m$ is the top degree homogeneous component of $k$. It can
be checked that this half-density is globally defined as an
element of the space $S^{m}(G_{{\mathcal F}_N},{\mathcal
L}(\pi^*E)\otimes |T{\mathcal G}_N|^{1/2})$.
\end{ex}

\subsection{Transverse bicharacteristic flow}\label{s:flow}
For any operator $P\in\Psi^m(M,E)$, let $\sigma_P$ denote the
transversal principal symbol of $P$, which is the restriction of
its principal symbol to $\tilde{N}^*{\mathcal F}$. We say that $P$
is transversally elliptic, if $\sigma_P(\nu)$ is invertible for
any $\nu\in\tilde{N}^*{\mathcal F}$.

Consider a transversally elliptic operator $A\in\Psi^2(M,E)$ which
has the scalar principal symbol and the holonomy invariant
transverse principal symbol. Here the holonomy invariance of the
transversal principal symbol $\sigma_{A}\in
C^{\infty}(\tilde{N}^*{\mathcal F})$ means that it is constant
along the leaves of the foliation $\cF_N$:
\[
\sigma_{A}(dh^*_{\gamma}(\nu))=\sigma_{A}(\nu),\quad \gamma\in G,
\quad \nu\in N^*_{r(\gamma)}\cF.
\]

Let $a_2\in S^2(\tilde{T}^*M)$ be the principal symbol of $A$.
(Here $\tilde{T}^*M= T^*M\setminus 0$). Take any scalar elliptic
symbol ${\tilde p}\in S^1(\tilde{T}^*M)$, which is equal to
$\sqrt{a_2}$ in some conic neighborhood of $\tilde{N}^*{\mathcal
F}$. Denote by $X_{\tilde p}$ the Hamiltonian vector field of
$\tilde p$ on $T^*M$. Since $N^*\cF$ is a coisotropic submanifold
in $T^*M$ and $T\cF_N$ is the symplectic orthogonal complement of
$T(N^*\cF)$, one can show that $X_{\tilde p}$ is tangent to
$\tilde{N}^*{\mathcal F}$, and its restriction to
$\tilde{N}^*{\mathcal F}$ (denoted also by $X_{\tilde p}$) is an
infinitesimal transformation of the foliation $\cF_N$, i.e. for
any vector field $X$ on $\tilde{N}^*{\mathcal F}$, tangent to
$\cF_N$, the commutator $[X_{\tilde p},X]$ is tangent to $\cF_N$.
It follows that the Hamiltonian flow $\tilde{f}_t$ of $\tilde{p}$
preserves $\tilde{N}^*{\mathcal F}$, and its restriction to
$N^*{\mathcal F}$ (denoted by $f_t$) preserves the foliation
$\cF_N$, that is, takes any leaf of $\cF_N$ to a leaf.

Let $\tau=T N^*\cF/T\cF_N$ be the normal space to the foliation
$\cF_N$ and $\pi_{tr}:T N^*\cF\to \tau$ the natural projection.
For any $(\gamma,\nu)\in G_{{\mathcal F}_N}$, let
$dH_{(\gamma,\nu)}: \tau_{dh^*_\gamma(\nu)} \to \tau_\nu $ be the
corresponding linear holonomy map. The differential of the map
$(s_N,r_N):G_{{\mathcal F}_N}\to N^*\cF\times N^*\cF$ at a point
$(\gamma,\nu)\in G_{{\mathcal F}_N}$ defines an inclusion of
$T_{(\gamma,\nu)} G_{{\mathcal F}_N}$ into $T_{dh^*_\gamma(\nu)}
N^*\cF \times T_\nu N^*\cF$, and its image consists of all
$(X,Y)\in T_{dh^*_\gamma(\nu)} N^*\cF \times T_\nu N^*\cF$ such
that
\begin{equation}\label{e:inf}
\pi_{tr}(Y) = dH_{(\gamma,\nu)}(\pi_{tr}(X)).
\end{equation}
Since $X_{\tilde p}$ is an infinitesimal transformation of the
foliation $\cF_N$, one can see that, for any $(\gamma,\nu)\in
G_{{\mathcal F}_N}$, the pair $(X_{\tilde p}(dh^*_\gamma(\nu)),
X_{\tilde p}(\nu)) \in T_{dh^*_\gamma(\nu)} N^*\cF \times T_\nu
N^*\cF$ satisfies (\ref{e:inf}). Therefore, there exists a unique
vector field $\cH_p$ on $G_{{\mathcal F}_N}$ such that
$ds_N(\cH_p)=X_{\tilde p}$ and $dr_N(\cH_p)=X_{\tilde p}$. Let
$F_t$ be the flow on $G_{{\mathcal F}_N}$ determined by the vector
field $\cH_p$. It is easy to see that $s_N\circ F_t=f_t\circ s_N$,
$r_N\circ F_t=f_t\circ r_N$ and the flow $F_t$ preserves the
foliation $\cG_N$.

\begin{defn}\label{d:flow} Let $P=\sqrt{A}$ be an (unbounded)
linear operator in $L^2(M,E)$, where $A\in\Psi^2(M,E)$ is an
essentially self-adjoint, transversally elliptic operator, which
has the scalar principal symbol and the holonomy invariant
transverse principal symbol. {\bf The transversal bicharacteristic
flow} of $P$ is the one-parameter group $F^*_t$ of automorphisms
of the involutive algebra $C^{\infty}_{prop}(G_{{\mathcal
F}_N},|T{\mathcal G}_N|^{1/2})$ induced by the flow $F_t$ on
$G_{{\mathcal F}_N}$.
\end{defn}

\begin{rem}
It is easy to see that the definition of transversal
bicharacteristic flow is independent of a choice of the elliptic
extension $\tilde{p}$.
\end{rem}

\begin{ex}
Consider a foliated coordinate chart $\varkappa : U\subset
M\stackrel{\sim}{\longrightarrow} I^{n}$ on $M$ with coordinates
$(x,y)\in I^p\times I^q$. Let $\tilde{p}$ be a positive, smooth
homogeneous of degree $1$ function on $I^n\times
(\RR^n\setminus\{0\})$ (a scalar elliptic principal symbol) such
that the corresponding transversal principal symbol $\sigma_P$ is
holonomy invariant. This means
\[
\tilde{p}(x,y,0,\eta)=p(y,\eta), \quad x\in I^p, \quad y\in I^q,
\quad \eta \in \RR^q
\]
with some function $p$. The Hamiltonian vector field $X_{\tilde
p}$ on $I^n\times \RR^n$ is given by
\[
X_{\tilde p} = \frac{\partial \tilde{p}}{\partial \xi}
\frac{\partial }{\partial x} - \frac{\partial \tilde{p}}{\partial
x} \frac{\partial }{\partial \xi} + \frac{\partial
\tilde{p}}{\partial \eta} \frac{\partial }{\partial y} -
\frac{\partial \tilde{p}}{\partial y} \frac{\partial }{\partial
\eta},
\]
and its restriction to $N^*\cF\left|_U\right.\cong I^p\times
I^q\times \RR^q$ is given by
\begin{multline*}
X_{\tilde p}(x,y,\eta) = \frac{\partial \tilde{p}}{\partial
\xi}(x,y,0,\eta) \frac{\partial }{\partial x} + \frac{\partial
p}{\partial \eta}(y,\eta) \frac{\partial }{\partial y} -
\frac{\partial p}{\partial y}(y,\eta) \frac{\partial }{\partial
\eta},\\ (x,y,\eta) \in I^p\times I^q\times \RR^q.
\end{multline*}
The fact that $X_{\tilde p}$ is an infinitesimal transformation of
the foliation $\cF_N$ means that its transverse part
\[
\frac{\partial p}{\partial \eta}(y,\eta) \frac{\partial }{\partial
y} - \frac{\partial p}{\partial y}(y,\eta) \frac{\partial
}{\partial \eta}
\]
is independent of $x$. The corresponding vector field $\cH_p$ on
$G_{{\mathcal F}_N}$ is given by
\begin{multline*}
\cH_p(x,x',y,\eta) = \frac{\partial \tilde{p}}{\partial
\xi}(x,y,0,\eta) \frac{\partial }{\partial x} + \frac{\partial
\tilde{p}}{\partial \xi}(x',y,0,\eta) \frac{\partial }{\partial
x'}\\ + \frac{\partial p}{\partial \eta}(y,\eta) \frac{\partial
}{\partial y} - \frac{\partial p}{\partial y}(y,\eta)
\frac{\partial }{\partial \eta},\quad (x,x',y,\eta) \in I^p\times
I^p\times I^q\times \RR^q.
\end{multline*}
Finally, the transversal bicharacteristic flow is given by the
action of the flow $F_t$ determined by the vector field $\cH_p$ on
the space of half-densities of the form $k_m(x,x',y,\eta)\,
|dx|^{1/2}\,|dx'|^{1/2}$.
\end{ex}

\begin{rem}
The construction of the transversal bicharacteristic flow provides
an example of what can be called noncommutative symplectic (or,
maybe, better, Poisson) reduction. Here symplectic reduction means
the following procedure \cite[Chapter III, Section 14]{LM87} (see
also \cite{Li75,Li77}).

Let $(X,\omega)$ be a symplectic manifold, and $Y$ a submanifold
of $X$ such that the $2$-form $\omega_Y$ induced by $\omega$ on
$Y$ is of constant rank. Let ${\mathcal F}_Y$ be the
characteristic foliation of $Y$ relative to $\omega_Y$. If the
foliation ${\mathcal F}_Y$ is simple, that is, it is given by the
fibers of a surjective submersion $p$ of $Y$ to a smooth manifold
$B$, then $B$ has a unique symplectic form $\omega_B$ such that
$p^*\omega_B=\omega_Y$. The symplectic manifold $(B, \omega_B)$ is
said to be the reduced symplectic manifold associated with $Y$. In
a particular case when the submanifold $Y$ is the preimage of a
point under the momentum map associated with the Hamiltonian
action of a Lie group, the symplectic reduction associated with
$Y$ is the Mardsen-Weinstein symplectic reduction \cite{MW}.

Moreover (see, for instance, \cite[Chapter III, Theorem
14.6]{LM87}), if $Y$ is invariant under the Hamiltonian flow of a
Hamiltonian $H\in C^\infty(X)$ (this is equivalent to the fact
that $(dH)\left|_Y\right.$ is constant along the leaves of the
characteristic foliation ${\mathcal F}_Y$), there exists a unique
function $\hat{H}\in C^\infty(B)$, called the reduced Hamiltonian,
such that $H\left|_Y\right. = \hat{H}\circ p$. Furthermore, the
map $p$ projects the restriction of the Hamiltonian flow of $H$ to
$Y$ to the reduced Hamiltonian flow on $B$ defined by the reduced
Hamiltonian $\hat{H}$.

Now let $(M,{\mathcal F})$ be a smooth foliated manifold. Consider
the symplectic reduction associated with the coisotropic
submanifold $Y=N^*{\mathcal F}$ in the symplectic manifold
$X=T^*M$. The corresponding characteristic foliation ${\mathcal
F}_Y$ is the linearized foliation ${\mathcal F}_N$. In general,
the leaf space $N^*{\cF}/{\mathcal F}_N$ is not a smooth manifold.
Following ideas of the noncommutative geometry in the sense of A.
Connes, one can treat the algebra $C^{\infty}_{prop}(G_{{\mathcal
F}_N},|T{\mathcal G}_N|^{1/2})$ as a noncommutative analogue of an
algebra of smooth functions on $N^*{\cF}/{\mathcal F}_N$. The
symplectic reduction procedure is applied to the Hamiltonian flow
$\tilde{f}_t$ of a function $\tilde p$ satisfying the assumptions
given in the beginning of this section, yielding the transversal
bicharacteristic flow $F^*_t$ as the corresponding reduced
Hamiltonian flow on $N^*{\cF}/{\mathcal F}_N$. Following the ideas
of \cite{Block-Ge, Xu}, one can interpret the algebra
$C^{\infty}_{prop}(G_{{\mathcal F}_N},|T{\mathcal G}_N|^{1/2})$ as
a noncommutative Poisson manifold and the flow $F^*_t$ as a
noncommutative Hamiltonian flow.
\end{rem}

\begin{ex}
Let $(M,{\mathcal F})$ be a compact Riemannian foliated manifold
equipped with a bundle-like metric $g_M$. Let $F=T{\mathcal F}$ be
the tangent bundle to ${\mathcal F}$, $H$ the orthogonal
complement to $F$, and $g_H$ the restriction of $g_M$ to $H$. By
definition, a Riemannian metric $g_M$ on $M$ is called
bundle-like, if it satisfies one of the following equivalent
conditions (see, for instance, \cite{Molino,Re}):
\medskip\par
1. For any continuous leafwise path $\gamma$ from $x$ to $y$, the
corresponding linear holonomy map $dh_\gamma : T_xM/T_x\cF\to
T_yM/T_y\cF$ is an isometry with respect to the Riemannian
structures on $T_xM/T_x\cF$ and $T_yM/T_y\cF$ induced by the
metric $g_M$;
\medskip\par
2. If $g_H$ is written as
$g_H=\sum_{\alpha\beta}g_{\alpha\beta}(x,y)\theta^\alpha\theta^\beta$
in some foliated chart with coordinates $(x,y)\in I^p\times I^q$,
where $\theta^\alpha\in H^*$ is the (unique) lift of $dy^\alpha$
under the projection $I^p\times I^q\to I^q$, then
$g_{\alpha\beta}$ is independent of $x$,
$g_{\alpha\beta}(x,y)=g_{\alpha\beta}(y)$.
\medskip\par
The decomposition $F\oplus H=TM$ induces a bigrading on $\bigwedge
T^{*}M$: $$ \bigwedge\nolimits^k
T^{*}M=\bigoplus_{i=0}^{k}\bigwedge\nolimits^{i,k-i}T^{*}M, $$
where $\bigwedge^{i,j}T^{*}M=\bigwedge^{i}F^{*}\otimes
\bigwedge^{j}H^{*}$. In this bigrading, the de Rham differential
$d$ can be written as $$ d=d_F+d_H+\theta, $$ where $d_F$ and
$d_H$ are first order differential operators (the tangential de
Rham differential and the transversal de Rham differential
accordingly), and $\theta$ is a zero order differential operator.

The transverse signature operator is a first order differential
operator in  $C^{\infty}(M,\bigwedge H^{*})$ given by $$ D_H=d_H +
d^*_H, $$ and the transversal Laplacian is a second order
transversally elliptic differential operator in
$C^{\infty}(M,\bigwedge H^{*})$ given by $$ \Delta_H=D_H^2. $$ The
principal symbol $\sigma(\Delta_H)$ of $\Delta_H$ is given by $$
\sigma(\Delta_H)(x,\xi)=g_H(\xi,\xi)I_x,\quad (x,\xi)\in
\tilde{T}^*M, $$ and holonomy invariance of the transversal
principal symbol is equivalent to the assumption on the metric
$g_M$ to be bundle-like.

Take any function $p_2\in C^{\infty}(T^*M)$, which coincides with
$\sqrt{\sigma(\Delta_H)}$ in some conical neighborhood of
$N^*{\mathcal F}$. The restriction of the Hamiltonian flow of
$p_2$ to $N^*{\mathcal F}$ coincides with the restriction $G_t$ of
the geodesic flow $g_t$ of the Riemannian metric $g_M$ to
$N^*{\mathcal F}$, which is the transversal bicharacteristic flow
of the operator $\langle D_H\rangle =\sqrt{\Delta_H+I}$.

Finally, if ${\mathcal F}$ is given by the fibers of a Riemannian
submersion $f:M\rightarrow B$, then there is a natural isomorphism
$N^*_m{\mathcal F}\rightarrow T^*_{f(m)}B$, and, under this
isomorphism, the transversal geodesic flow $G_t$ on $N^*{\mathcal
F}$ corresponds to the geodesic flow $T^*B$ (see, for instance,
\cite{ONeil,Re}).
\end{ex}

\subsection{Egorov's theorem}\label{s:Egorov}
Let $D\in\Psi^1(M,E)$ be a formally self-adjoint, transversally
elliptic operator such that $D^2$ has the scalar principal symbol
and the holonomy invariant transverse principal symbol. By
\cite{noncom}, the operator $D$ is essentially self-adjoint with
initial domain $C^\infty(M,E)$. Define an unbounded linear
operator $\langle D\rangle$ in the space $L^2(M,E)$ as
\[
\langle D\rangle=(D^2+I)^{1/2}.
\]
By the spectral theorem, the operator $\langle D\rangle $ is
well-defined as a positive, self-adjoint operator in $L^2(M,E)$.
The operator $\langle D\rangle^2\in \Psi^2(M,E)$ is a bounded
operator from $H^2(M,E)$ to $L^2(M,E)$. Hence, by interpolation,
$\langle D\rangle$ defines a bounded operator from $H^1(M,E)$ to
$L^2(M,E)$ and $H^1(M,E)$ is contained in the domain of $\langle
D\rangle$ in $L^2(M,E)$.

By the spectral theorem, the operator $\langle
D\rangle^s=(D^2+I)^{s/2}$ is a well-defined positive self-adjoint
operator in $\cH=L^2(M,E)$ for any $s\in\RR$, which is unbounded
if $s>0$. For any $s\geq 0$, denote by $\cH^s$ the domain of
$\langle D\rangle^s$, and, for $s<0$, $\cH^s=(\cH^{-s})^*$. Put
also $\cH^{\infty}=\bigcap_{s\geq 0}\cH^s, \quad
\cH^{-\infty}=(\cH^{\infty})^*$. It is clear that $H^s(M,E)
\subset \cH^s$ for any $s\geqslant 0$ and $\cH^s \subset H^s(M,E)$
for any $s<0$.  In particular, $C^\infty(M,E) \subset \cH^s$ for
any $s$.

We say that a bounded operator $A$ in $\cH^{\infty}$ belongs to
$\cL(\cH^{-\infty},\cH^{\infty})$ (resp.
$\cK(\cH^{-\infty},\cH^{\infty})$), if, for any $s$ and $r$, it
extends to a bounded (resp. compact) operator from $\cH^s$ to
$\cH^r$, or, equivalently, the operator $\langle
D\rangle^rA\langle D\rangle^{-s}$ extends to a bounded (resp.
compact) operator in $L^2(M,E)$. It is easy to see that
$\cL(\cH^{-\infty},\cH^{\infty})$ is a involutive subalgebra in
$\cL(\cH)$ and $\cK(\cH^{-\infty},\cH^{\infty})$ is its ideal. We
also introduce the class $\cL^1(\cH^{-\infty},\cH^{\infty})$,
which consists of all operators from
$\cK(\cH^{-\infty},\cH^{\infty})$ such that, for any $s$ and $r$,
the operator $\langle D\rangle^rA\langle D\rangle^{-s}$  is a
trace class operator in $L^2(M,E)$. It should be noted that any
operator $K$ with the smooth kernel belongs to
$\cL^1(\cH^{-\infty},\cH^{\infty})$.

As an operator acting on half-densities, any operator
$P\in\Psi^m(M)$ has the subprincipal symbol which is the
well-defined homogeneous of degree $m-1$ smooth function on
$T^*M\setminus 0$ given in local coordinates by the formula
\begin{equation}\label{e:subprincipal1}
p_{sub}=p_{m-1}-\frac{1}{2i}\sum_{j=1}^n\frac{\partial^2p_m}{\partial
x_j\partial \xi_j},
\end{equation}
where $p_{m-1}$ and $p_m$ are the homogeneous components of the
complete symbol of $P$ of degree $m-1$ and $m$ respectively.
Observe that $p_{sub}=0$ if $P$ is a real, self-adjoint,
differential operator of even order. In particular, this holds for
the transversal Laplacian $\Delta_H$ on functions.

By the spectral theorem, the operator $\langle D\rangle$ defines a
strongly continuous group $e^{it\langle D\rangle}$ of bounded
operators in $L^2(M,E)$. Consider a one-parameter group $\Phi_t$
of $\ast$-automorphisms of the algebra ${\mathcal L}(L^2(M,E))$
defined by
\begin{displaymath}
\Phi_t(T)=e^{i t\langle D\rangle}Te^{-i t\langle D\rangle}, \quad
T\in {\mathcal L}(L^2(M,E)).
\end{displaymath}

The main result of the paper is the following theorem.

\begin{thm}
\label{Egorov} Let $D\in\Psi^1(M,E)$ be a formally self-adjoint,
transversally elliptic operator such that $D^2$ has the scalar
principal symbol and the holonomy invariant transverse principal
symbol.
\medskip\par
(1) For any $K\in \Psi^{m,-\infty}(M,{\mathcal F},E)$, there
exists an operator $K(t)\in\Psi^{m,-\infty}(M,{\mathcal F},E)$
such that $\Phi_t(K)-K(t), t\in \RR,$ is a smooth family of
operators of class $\cL^1(\cH^{-\infty},\cH^{\infty})$.
\medskip\par
(2) If, in addition, $E$ is the trivial line bundle, and the
subprincipal symbol of $D^2$ is zero, then, for any $K\in
\Psi^{m,-\infty}(M,{\mathcal F})$ with the principal symbol $k\in
S^{m}(G_{{\mathcal F}_N},|T{\mathcal G}_N|^{1/2})$, the principal
symbol $k(t)\in S^m(G_{{\mathcal F}_N},|T{\mathcal G}_N|^{1/2})$
of the operator $K(t)$ is given by $k(t)=F^*_t(k)$, where $F^*_t$
is the transverse bicharacteristic flow of the operator $\langle
D\rangle $.
\end{thm}

\begin{rem}
Theorem~\ref{Egorov} implies Egorov's theorem for elliptic
operators on compact Riemannian orbifolds. An $m$-dimensional
orbifold $M$ is a Hausdorff, second countable topological space,
which is locally diffeomorphic to the quotient of $\RR^m$ by a
finite group of diffeomorphisms $\Gamma$. The notion of orbifold
was first introduced by Satake in \cite{Satake}, where a different
name, $V$-manifold, was used. We refer the reader to
\cite{Satake,Ka78,Chen-Ruan} for expositions of orbifold theory.
It is well-known (see, for instance, \cite{Ka81}) that any
orbifold $M$ is diffeomorphic to the orbifold of $G$ orbits of an
action of a compact Lie group $G$ on a compact manifold $P$ where
the action has finite isotropy groups (actually, one can take $P$
to be the orthogonal frame bundle of $M$ and $G=O(m)$). The orbits
of this action are the leaves of a foliation $\cF$ on $P$. We will
use a natural isomorphism of the space $C^\infty(M)$ with the
space $C^\infty(P)^G$ of $G$ invariant functions on $P$. A
pseudodifferential operator $A$ in $C^\infty(P)$ can be defined as
an operator acting on $C^\infty(P)^G$ which is the restriction of
a $G$ equivariant pseudodifferential operator $\tilde{A}$ in
$C^\infty(M)$. The operator $A$ is elliptic iff the corresponding
operator $\tilde{A}$ is transversally elliptic with respect to the
foliation $\cF$. The orthogonal projection $\Pi$ on the space of
$G$-invariant functions in $C^\infty(P)$ is a transversal
pseudodifferential operator of class $\Psi^{0,-\infty}(P,\cF)$. It
follows that a pseudodifferential operator $A$ in $C^\infty(M)$
coincides with the restriction of the operator $\Pi\tilde{A}\Pi\in
\Psi^{0,-\infty}(P,\cF)$ to $C^\infty(P)^G$.

Fix Riemannian metrics $g_M$ on $M$ and $g_P$ on $P$ such that the
quotient map $P\to M$ is a Riemannian submersion. So $g_P$ is a
bundle-like metric on the foliated manifold $(P,\cF).$ One can
show that the associated transverse Laplacian $\Delta_H$ is
$G$-invariant and the Laplacian $\Delta_M$ on $M$ coincides with
the restriction of $\Delta_H$ to $C^\infty(P)^G$. Therefore, we
have
\begin{align*}
e^{it(\Delta_M +I)^{1/2}} A e^{-it (\Delta_M +I)^{1/2}} & = \Pi
e^{it (\Delta_H +I)^{1/2}}\tilde{A}e^{-it(\Delta_H +I)^{1/2}} \Pi
\\ & = \Pi e^{it (\Delta_H +I)^{1/2}}(\Pi\tilde{A}\Pi
)e^{-it(\Delta_H +I)^{1/2}} \Pi.
\end{align*}
By Theorem~\ref{Egorov}, it follows that the operator
$e^{it(\Delta_M +I)^{1/2}} A e^{-it (\Delta_M +I)^{1/2}}$ is a
pseudodifferential operator on $M$ and one can describe its
principal symbol as in the classical Egorov's theorem. The details
will be given elsewhere.
\end{rem}

\subsection{Noncommutative geodesic flow on foliated manifolds}
\label{ncg}

As stated in \cite{noncom}, any operator $D$, satisfying the
assumptions of Section~\ref{s:Egorov}, defines a spectral triple
in the sense of Connes' noncommutative geometry. In this setting,
Theorem~\ref{Egorov} has a natural interpretation in terms of the
corresponding noncommutative geodesic flow. First, we recall
general definitions \cite{Co-M,Sp-view}.

Let $({\mathcal A},{\mathcal H},D)$ be a spectral triple
\cite{Sp-view}. Here
\begin{enumerate}
\item ${\mathcal A}$ is an involutive algebra;

\item ${\mathcal H}$ is a Hilbert space equipped with a
$\ast$-representation of the algebra ${\mathcal A}$ (we will
identify an element $a\in \cA$ with the corresponding operator in
$\mathcal H$);

\item $D$ is an (unbounded) self-adjoint operator in ${\mathcal H}$ such
that

\begin{description}
\item[(a)] for any $a\in {\mathcal A}$, the operator $a(D-i)^{-1}$ is a
compact operator in ${\mathcal H}$;

\item[(b)] $D$ almost commutes with any $a\in {\mathcal A}$ in the sense
that $[D,a]$ is bounded in $\cH$.
\end{description}
\end{enumerate}

As above, let $\langle D\rangle=(D^2+I)^{1/2}$. By $\delta$, we
denote the (unbounded) derivative on ${\mathcal L}(\cH)$ given by
\begin{equation}\label{e:delta}
\delta(T)=[\langle D\rangle,T], \quad T\in \Dom\delta\subset
{\mathcal L}(\cH).
\end{equation}
Let ${\rm OP}^{\alpha}$ be the space of operators in ${\mathcal
H}$ of order ${\alpha}$, that means that $P\in {\rm OP}^{\alpha}$
iff $P\langle D\rangle^{-\alpha}\in \bigcap_n {\rm
Dom}\;\delta^n$. In particular, ${\rm OP}^{0}=\bigcap_n {\rm
Dom}\;\delta^n$. Denote by ${\rm OP}_0^{0}$ the space of all
operators $P\in {\rm OP}^{0}$ such that $\langle D\rangle^{-1}P$
and $P\langle D\rangle^{-1}$ are compact operators in $\cH$. We
also say that $P\in {\rm OP}_0^{\alpha}$ if $P\langle
D\rangle^{-\alpha}$ and $\langle D\rangle^{-\alpha}P$ are in ${\rm
OP}_0^{0}$. It is easy to see that ${\rm
OP}_0^{-\infty}=\bigcap_{\alpha}{\rm OP}_0^{\alpha}$ coincides
with ${\mathcal K}({\mathcal H}^{-\infty}, {\mathcal H}^\infty)$.

We will assume that $({\mathcal A},{\mathcal H},D)$ is smooth.
This means that, for any $a\in {\mathcal A}$, the bounded
operators $a$ and $[D,a]$ in ${\mathcal H}$ belong to ${\rm
OP}^{0}$. Let ${\mathcal B}$ be the algebra of bounded operators
in ${\mathcal H}$ generated by the set of all operators of the
form $\delta^n(a)$ with $a\in {\mathcal A}$ and $n\in {\NN}$.
Furthermore, we assume that the algebra ${\mathcal B}$ is
contained in ${\rm OP}_0^{0}$. In particular, this implies that
$({\mathcal B},{\mathcal H},D)$ is a spectral triple in the above
sense.

In \cite{Co-M,Sp-view}, the definition of the algebra
$\Psi^*({\mathcal A})$ of pseudodifferential operators was given
for a unital algebra $\cA$. In the case under consideration, the
algebra ${\mathcal A}$ is non-unital, that, roughly speaking,
means that the associated geometric space is noncompact.
Therefore, we must take into account behavior of
pseudodifferential operators at "infinity". Next we define an
algebra $\Psi^*_0({\mathcal A})$, which can considered as an
analogue of the algebra of pseudodifferential operators on a
noncompact Riemannian manifold, whose symbols and all its
derivatives of any order vanish at infinity. In particular, the
assumptions on the spectral triple made above mean that the
algebra $\cA$ consists of smooth "functions", vanishing at
"infinity" with all its derivatives of any order.

Define $\Psi^*_0({\mathcal A})$ as the set of (unbounded)
operators in $\cH$, which admit an asymptotic expansion:
\begin{equation}
\label{psif} P\sim \sum_{j=0}^{+\infty}b_{q-j}\langle
D\rangle^{q-j}, \quad b_{q-j}\in {\mathcal B},
\end{equation}
that means that, for any $N$,
\[
P - \left(b_q\langle D\rangle^q + b_{q-1}\langle
D\rangle^{q-1}+\ldots+b_{-N}\langle D\rangle^{-N}\right)\in {\rm
OP}_0^{-N-1}.
\]
By an easy modification of the proof of Theorem B.1 in
\cite[Appendix B]{Co-M}, one can prove that $\Psi^*_0({\mathcal
A})$ is an algebra. Let ${\mathcal C}_0$ be the algebra ${\mathcal
C}_0={\rm OP}_0^0\bigcap \Psi^*_0({\mathcal A})$, and
$\bar{\mathcal C}_0$ the closure of ${\mathcal C}_0$ in ${\mathcal
L}(\cH)$.

For any $T\in \cL(\cH)$, define
\begin{equation}\label{e:alphat}
\alpha_t(T)=e^{it\langle D\rangle}Te^{-it\langle D\rangle}, \quad
t\in {\RR}.
\end{equation}
As usual, ${\mathcal K}$ denotes the ideal of compact operators in
${\mathcal H}$. The following definitions are motivated by the
work of Connes \cite{Sp-view}.

\begin{defn}
Under the current assumptions on a spectral triple $({\mathcal A},
{\mathcal H}, D)$, {\bf the unitary cotangent bundle}
$S^*{\mathcal A}$ is defined as the quotient of the $C^*$-algebra
generated by all $\alpha_t(\bar{\mathcal C}_0), t\in\RR$ and
$\mathcal K$ by ${\mathcal K}$.
\end{defn}

\begin{defn}
Under the current assumptions on a spectral triple $({\mathcal A},
{\mathcal H}, D)$, {\bf the noncommutative geodesic flow} is the
one-pa\-ra\-me\-ter group $\alpha_t$ of automorphisms of the
algebra $S^*{\mathcal A}$ defined by (\ref{e:alphat}).
\end{defn}

We consider spectral triples $({\mathcal A},{\mathcal H},D)$
associated with a compact foliated Riemannian manifold
$(M,{\mathcal F})$ \cite{noncom}:
\begin{enumerate}
\item The involutive
algebra ${\mathcal A}$ is the algebra $\cinf_c(G,|T\cG|^{1/2})$;
\item The Hilbert space ${\mathcal H}$
is the space $L^2(M,E)$ of $L^2$-sections of a holonomy
equivariant Hermitian vector bundle $E$, on which an element $k$
of the algebra ${\mathcal A}$ is represented via the
$\ast$-representation $R_E$ (see below for a definition);
\item The operator $D$ is a first order self-adjoint
transversally elliptic operator with the holonomy invariant
transversal principal symbol such that the operator $D^2$ has the
scalar principal symbol.
\end{enumerate}

We recall briefly the definitions of the structure of involutive
algebra on $\cA$ and of the representation $R_E$. Let $\alpha\in
C^{\infty}(M,|T{\cF}|^{1/2})$ be a strictly positive, smooth,
leafwise half-density. One can lift $\alpha$ to a strictly
positive, leafwise half-density $\nu^x=s^*\alpha\in
C^{\infty}(G^x,|TG^x|^{1/2})$ via the covering map $s:G^x\to L_x$
($L_x$ is the leaf through a point $x\in M$). In the presence of
$\nu$, the space $\cA=C^{\infty}_c(G,|T{\cG}|^{1/2})$ is naturally
identified with $C^{\infty}_c(G)$. We also assume, for simplicity,
that there exists a holonomy invariant, smooth, transverse
half-density $\Lambda\in C^{\infty}(M,|TM/T\cF|^{1/2})$. Recall
that the holonomy invariance of $\Lambda$ means that
$dh^*_{\gamma}(\Lambda(y)) = \Lambda(x)$ for any $\gamma\in G,
s(\gamma)=x, r(\gamma)=y$, where the map $dh^*_{\gamma}:
|T_yM/T_y{\mathcal F}|^{1/2}\to |T_xM/T_x{\mathcal F}|^{1/2}$ is
induced by the corresponding linear holonomy map.

The multiplication and the involution in $\cA$ are given by the
formulas
\begin{equation}\label{e:a}
\begin{aligned}
(k_1\ast k_2)(\gamma) &
=\int_{G^x}k_1(\gamma'^{-1}\gamma)k_2(\gamma')\,d\nu^x(\gamma'),
\quad \gamma\in G^x,\\ k^*(\gamma) &
=\overline{k(\gamma^{-1})},\quad \gamma\in G,
\end{aligned}
\end{equation}
where $k,k_1, k_2\in \cA$.

An Hermitian vector bundle $E$ on $M$ is holonomy equivariant, if
it is equipped with an isometric action $$
T(\gamma):E_x\rightarrow E_y, \quad \gamma\in G,
\gamma:x\rightarrow y, $$ of $G$ in fibers of $E$. Using the fixed
half-densities $\alpha$ and $\Lambda$, one can identify elements
of $L^2(M,E)$ with square integrable sections of the bundle $E$.
Then, for any $u\in L^2(M,E)$, the section $R_E(k)u\in L^{2}(M,E)$
is defined by the formula $$ R_E(k)u(x) =\int_{ G^{x}}
k(\gamma)\,T(\gamma)[u(s(\gamma))]\,d\nu^x(\gamma),\quad x\in M.
$$

It was stated in \cite{noncom} that the spectral triple
$({\mathcal A},{\mathcal H},D)$ associated with a compact foliated
Riemannian manifold is smooth. Recall that this means that, for
any $a\in {\mathcal A}$, $a$ and $[D,a]$ belong to ${\rm
OP}^{0}=\bigcap_n {\rm Dom}\;\delta^n$. There is a gap in the
proof of this fact given in \cite{noncom}. In this paper, we give
a correct proof (cf. Theorem~\ref{p:smooth} below). In
Theorem~\ref{p:smooth}, we also prove that, in the case in
question, the algebra ${\mathcal B}$ mentioned above is contained
in ${\rm OP}_0^{0}$.

For any $\nu\in \tilde{N}^*\cF$, there is a natural
$\ast$-representation $R_\nu$ of the algebra $S^{0}(G_{{\mathcal
F}_N}, |T{\mathcal G}_N|^{1/2})$ in
$L^2(G^\nu_{{\cF}_N},s_N^*(\pi^*E))$. For its definition, we will
use the strictly positive, leafwise half-density $\mu^\nu\in
C^{\infty}(G^\nu_{{\cF}_N}, |T{\mathcal G}_N|^{1/2})$ induced by
$\alpha$ and the corresponding isomorphism $S^{0}(G_{{\mathcal
F}_N}, |T{\mathcal G}_N|^{1/2})\cong S^{0}(G_{{\mathcal F}_N})$.
Since $E$ is a holonomy equivariant vector bundle, the bundle
$\pi^*E$ is also holonomy equivariant. The action of $G_{{\mathcal
F}_N}$ in fibers of $\pi^*E$,
\begin{displaymath}
\pi^*T(\gamma,\nu):(\pi^*E)_{dh^*_{\gamma}(\nu)}\rightarrow
(\pi^*E)_\nu,\quad (\gamma,\nu)\in G_{{\mathcal F}_N},
\end{displaymath}
is given by the formula  $\pi^*T(\gamma,\nu)=T(\gamma)$, where we
use the natural isomorphisms $(\pi^*E)_{dh^*_{\gamma}(\nu)}=E_x$
and $(\pi^*E)_\nu=E_y$. For any $k\in S^{0}(G_{{\mathcal F}_N})$
and $u\in L^2(G^\nu_{{\cF}_N}, s_N^*(\pi^*E))$, the section
$R_\nu(k)u\in L^2(G^\nu_{{\cF}_N}, s_N^*(\pi^*E))$ is given by the
formula
\begin{multline*}
R_\nu(k)u(\gamma,\nu)
=\int_{G^\nu_{{\cF}_N}}k((\gamma',\nu)^{-1}(\gamma,\nu))
\pi^*T(\gamma',\nu)[u(\gamma',\nu)]\,d\mu^\nu(\gamma',\nu),\\
(\gamma,\nu)\in G^\nu_{{\cF}_N}.
\end{multline*}

It follows from the direct integral decomposition
\[
L^2(G_{{\cF}_N},s_N^*(\pi^*E))=
\int_{N^*\cF}L^2(G^\nu_{{\cF}_N},s_N^*(\pi^*E))\,d\nu,
\]
that, for any $k\in S^{0}(G_{{\mathcal F}_N}, |T{\mathcal
G}_N|^{1/2})$, the continuous family
\[\{R_\nu(k)\in\cL(L^2(G^\nu_{{\cF}_N},s_N^*(\pi^*E))):\nu\in
\tilde{N}^*\cF\}\] defines a bounded operator in
$L^2(G_{{\cF}_N},s_N^*(\pi^*E))$. We will identify $k\in
S^{0}(G_{{\mathcal F}_N}, |T{\mathcal G}_N|^{1/2})$ with the
corresponding operator in $L^2(G_{{\cF}_N},s_N^*(\pi^*E))$ and
denote by $\bar{S}^{0}(G_{{\mathcal F}_N}, |T{\mathcal
G}_N|^{1/2})$ the closure of $S^{0}(G_{{\mathcal F}_N},
|T{\mathcal G}_N|^{1/2})$ in the uniform operator topology of
$\cL(L^2(G_{{\cF}_N},s_N^*(\pi^*E)))$. The transversal
bicharacteristic flow $F^*_t$ of the operator $\langle D\rangle$
extends by continuity to a strongly continuous one-parameter group
of automorphisms of $\bar{S}^0(G_{{\mathcal
F}_N},|T{\cG}_N|^{1/2})$.

The following theorem gives a description of the associated
noncommutative geodesic flow in the scalar case.

\begin{thm}
\label{noncom:flow} Let $({\mathcal A},{\mathcal H},D)$ be a
spectral triple associated with a compact foliated Riemannian
manifold $(M,{\mathcal F})$ as above with $E$, being the trivial
holonomy equivariant line bundle. Assume that the subprincipal
symbol of $D^2$ vanishes. There exists a surjective homomorphism
of involutive algebras $P : S^*{\mathcal A}\rightarrow
\bar{S}^0(G_{{\mathcal F}_N},|T{\cG}_N|^{1/2})$ such that the
following diagram commutes:
\begin{equation}\label{e:cd}
  \begin{CD}
S^*{\mathcal A} @>\alpha_t>> S^*{\mathcal A}\\ @VPVV @VVPV
\\ \bar{S}^0(G_{{\mathcal F}_N},|T{\cG}_N|^{1/2})@>F^*_t>>
\bar{S}^0(G_{{\mathcal F}_N},|T{\cG}_N|^{1/2})
  \end{CD}
\end{equation}
\end{thm}

\section{Proof of the main theorem}
\label{s:proof}
\subsection{The case of elliptic operator}
\label{elliptic:sect} Let $(M,{\mathcal F})$ be a compact foliated
manifold, $E$ a Hermitian vector bundle on $M$. In this section,
we will assume that $D\in \Psi^1(M,E)$ is a formally self-adjoint,
{\em elliptic} operator such that $D^2$ has the scalar principal
symbol and the holonomy invariant transverse principal symbol.
Then $P=\langle D\rangle\in \Psi^1(M,E)$ is a self-adjoint
elliptic operator with the positive, scalar principal symbol $p$
and the holonomy invariant transversal principal symbol. In this
case, the elliptic extension $\tilde p$ of $p$ introduced in
Section~\ref{s:flow} can be taken to be equal to $p$, $\tilde
p=p$. Therefore, if we denote by $X_{p}$ the Hamiltonian vector
field of $p$ on $T^*M$, then the vector field $\cH_p$ can be
described as a unique vector field on $G_{{\mathcal F}_N}$ such
that $ds_N(\cH_p)=X_{p}$ and $dr_N(\cH_p)=X_{p}$. Similarly, one
can define the transverse bicharacteristic flow $F^*_t$ of $P$ as
in Definition~\ref{d:flow}, using $p$ instead of $\tilde p$. The
following theorem is slightly stronger than Theorem~\ref{Egorov}.

\begin{thm}
\label{Egorov:elliptic} For any $K\in \Psi^{0,-\infty}(M,{\mathcal
F},E)$, the operator \[\Phi_t(K)=e^{itP}Ke^{-itP}\] is an operator
of class $\Psi^{0,-\infty}(M,{\mathcal F},E)$.

If $E$ is the trivial line bundle, and the subprincipal symbol of
$D^2$ vanishes, then, for any operator $K\in
\Psi^{0,-\infty}(M,{\mathcal F})$ with the principal symbol $k\in
S^{0}(G_{{\mathcal F}_N},|T{\mathcal G}_N|^{1/2})$, the operator
$\Phi_t(K)$ has the principal symbol $k(t)\in S^0(G_{{\mathcal
F}_N},|T{\mathcal G}_N|^{1/2})$ given by $k(t)=F^*_t(k)$.
\end{thm}

\begin{proof}
For the proof, we use theory of Fourier integral operators (see,
for instance, \cite{H4,Taylor,Treves2}). Recall that a Fourier
integral operator on $M$ is a linear operator $F:C^{\infty}(M)\to
{\cD}'(M)$, represented microlocally in the form
\begin{equation}\label{e:FIO}
Fu(x)=\int e^{\phi(x,y,\theta)}a(x,y,\theta)\,u(y)\,dy\,d\theta,
\end{equation}
where $x\in X\subset \RR^n, y\in Y\subset\RR^n, \theta\in
\RR^N\setminus 0$. Here $a(x,y,\theta)\in S^m(X\times
Y\times\RR^N)$ is an amplitude, $\phi$ is a non-degenerate phase
function.

Consider the smooth map from $X\times Y\times\RR^N$ to $T^*X\times
T^*Y$ given by
\[
(x,y,\theta)\mapsto (x,\phi_x(x,y,\theta),y,-\phi_y(x,y,\theta)).
\]
The image of the set
\[
\Sigma_\phi=\{(x,y,\theta)\in X\times Y\times\RR^N
:\phi_\theta(x,y,\theta)=0\}
\]
under this map turns out to be a homogeneous canonical relation
$\Lambda_\phi$ in $T^*X\times T^*Y$. (Recall that a closed conic
submanifold $C\in T^*(X\times Y)\setminus 0$ is called a
homogeneous canonical relation, if it is Lagrangian with respect
to the 2-form $\omega_X-\omega_Y$, where $\omega_X, \omega_Y$ are
the canonical symplectic forms in $T^*X, T^*Y$ accordingly.)

The Fourier integral operator $F$ given by the formula
(\ref{e:FIO}) is said to be associated with $\Lambda_\phi$. We
will write $F\in I^m(X\times Y, \Lambda_\phi)$, if $a\in
S^{m+n/2-N/2}(X\times Y\times\RR^N)$.

Operators from $\Psi^{m,-\infty}(M,{\mathcal F},E)$ can be
described as Fourier integral operators associated with the
immersed canonical relation $G'_{{\mathcal F}_N}$, which is the
image of $G_{{\mathcal F}_N}$ under the mapping $G_{{\mathcal
F}_N}\rightarrow T^*M\times T^*M: (\gamma,\nu)\mapsto
(r_N(\gamma,\nu), -s_N(\gamma,\nu))$ \cite{noncom}. Indeed,
consider an elementary operator $A:C^{\infty}_c(U,\left.
E\right|_U)\to C^{\infty}_c(U',\left. E\right|_{U'})$ given by the
formula (\ref{loc}) with $k \in S ^{m} (I ^{p} \times I^p\times
I^q\times {\RR}^{q}, {\cL}({\CC}^r))$. It can be represented in
the form (\ref{e:FIO}), if we take $X=U$ with coordinates $(x,y)$,
$Y=U'$ with coordinates $(x',y')$, $\theta=\eta $, $N=q$, a phase
function $\phi(x,y,x',y')=(y-y')\eta$ and an amplitude
$a=k(x,x',y,\eta)$. The associated homogeneous canonical relation
$\Lambda_\phi$ is the set of all
$(x,y,\xi,\eta,x',y',\xi',\eta')\in T^*U\times T^*U'$ such that
$y=y',\xi=\xi'=0,\eta=-\eta'$, that coincides with the
intersection of $G'_{{\mathcal F}_N}$ with $T^*U\times T^*U'$.
Moreover, we see that
\[
\Psi^{m,-\infty}(M,{\mathcal F},E)\subset I^{m-p/2}(M\times
M,G'_{{\mathcal F}_N};\cL(E)\otimes |T(M\times M)|^{1/2}).
\]
Since $G_{{\mathcal F}_N}$ is, in general, an immersed canonical
relation, it is necessary to be more precise in the definition of
the classes $I^{m}(M\times M,G'_{{\mathcal F}_N};\cL(E)\otimes
|T(M\times M)|^{1/2})$. This can be done by analogy with the
definition of the classes of longitudinal pseudodifferential
operators on a foliated manifold given in \cite{Co79} (see also
\cite{noncom} and the definition of classes
$\Psi^{0,-\infty}(M,\cF,E)$ given above).

Let $p$ be the principal symbol of $P$, and let $\Lambda_{p}(t),
t\in\RR,$ be the canonical relation in $T^*M \times {T}^*M$
defined as
\[
\Lambda_{p}(t)= \{((x,\xi),(y,\eta))\in T^*M\times T^*M:
(x,\xi)=f_{-t}(y,\eta)\},
\]
where $f_t$ is the Hamiltonian flow of $p$. It is well-known (cf.,
for instance, \cite{Taylor}) that $e^{itP}$ is a Fourier integral
operator associated with $\Lambda_p(t)$:
\[
e^{itP}\in I^0(M\times M, \Lambda_{{p}}(t);\cL(E)\otimes
|T(M\times M)|^{1/2}).
\]
By holonomy invariance of the transverse principal symbol of $P$,
it follows that $\Lambda_{p}(t)\circ G_{{\mathcal
F}_N}\circ\Lambda_{p}(-t)=G_{{\mathcal F}_N}$, and by the
composition theorem of Fourier integral operators (see, for
instance, \cite{H4}), we have $\Phi_t(K)=e^{itP}Ke^{-itP}\in
\Psi^{0,-\infty}(M,{\mathcal F},E)$.

Now assume, in addition, that $E$ is the trivial line bundle, the
subprincipal symbol of $D^2$ vanishes, and $K\in
\Psi^{0,-\infty}(M,{\mathcal F})$ with the principal symbol $k\in
S^{0}(G_{{\mathcal F}_N},|T{\mathcal G}_N|^{1/2})$. Denote by
$\cL_{\cH_p}$ the Lie derivative on $C^{\infty}(G_{{\mathcal
F}_N},|T{\mathcal G}_N|^{1/2})$ by the vector field $\cH_p$. So
the function $k(t)=F^*_t(k) \in S^{0}(G_{{\mathcal
F}_N},|T{\mathcal G}_N|^{1/2})$ is the solution of the equation $$
\frac{d k(t)}{d t}=\cL_{\cH_p} k(t),\quad t\in\RR, $$ with the
initial data $k(0)=k$. By \cite{GS79} (cf. also
\cite{Duist-Ho2,H4}), it follows that, for any $K_1\in
\Psi^{0,-\infty}(M,{\mathcal F})$, the operator $[P,K_1]$ belongs
to $\Psi^{0,-\infty}(M,{\mathcal F})$, and
\begin{displaymath}
\sigma([P, K_1])=\frac{1}{i}\cL_{\cH_p} \sigma(K_1).
\end{displaymath}
Consider any smooth family $\cK(t)\in \Psi^{0,-\infty}(M,\cF),
t\in\RR,$ of operators with the principal symbol $k(t)$. Then
\begin{align*}
\frac{d \cK(t)}{d t}&=i[P,\cK(t)]+R(t),\quad t\in\RR,\\
\cK(0)&=K+R_0,
\end{align*}
where $R(t)\in \Psi^{-1,-\infty}(M,\cF), t\in\RR,$ is a smooth
family of operators, and $R_0\in \Psi^{-1,-\infty}(M,\cF)$.

Using the fact that $\Phi_t(K)$ is the solution of the Cauchy
problem
\begin{align*}
\frac{d \Phi_t(K)}{d t}&=i[P,\Phi_t(K)],\quad t\in\RR,\\
\Phi_0(K)&=K,
\end{align*}
and the first part of the theorem, we get
\begin{displaymath}
\cK(t)-\Phi_t(K)=\int_0^t \Phi_{t-\tau}(R(\tau))\,d\tau +
\Phi_t(R_0)\in\Psi^{-1,-\infty}(M,\cF),
\end{displaymath}
and $\sigma(\Phi_t(K)) =\sigma(\cK(t)) =k(t)$.
\end{proof}

\subsection{The general case}
\label{general:sect} In this section, we will prove
Theorem~\ref{Egorov} in the general case. Thus, we assume that
$D\in\Psi^1(M,E)$ is a formally self-adjoint, transversally
elliptic operator such that $D^2$ has the scalar principal symbol
and the holonomy invariant transverse principal symbol.

\begin{defn}
An operator $A\in \Psi^{l}(M,E)$ is said to be of order $-\infty$
in some conic neighborhood of $N^{*}\cF$, if, in any regular
foliated chart with the coordinates $(x,y)\in I^p\times I^q$,
there exists $\varepsilon>0$ such that, for any multiindices
$\alpha$ and $\beta$ and for any natural $N$, its complete symbol
$a\in S^l(I^n\times \RR^n)$ satisfies the estimate with some
constant $C_{\alpha\beta N}>0$
\begin{multline*}
|\partial^{\alpha}_\xi\partial^{\beta}_xa(x,y,\xi,\eta)|<C_{\alpha\beta
N} (1+|\xi|+|\eta|)^{-N},\\ (x,y)\in I^p\times I^q, \quad
(\xi,\eta)\in \RR^p\times \RR^q, \quad |\xi|< \varepsilon |\eta|.
\end{multline*}
\end{defn}

The important fact, concerning to operators of order $-\infty$ in
some conic neighborhood of $N^{*}\cF$, is contained in the
following lemma \cite{noncom}:

\begin{lem}
\label{smooth} If $A\in \Psi^{l}(M,E)$ is of order $-\infty$ in
some conic neighborhood of $N^{*}\cF$ and $K\in
\Psi^{m,-\infty}(M,{\mathcal F},E)$, then $AK$ and $KA$ are in
$\Psi^{-\infty}(M,E)$.
\end{lem}

Denote by $\cL(\cD'(M,E),\cH^{\infty})$ (resp.
$\cL(\cH^{-\infty},C^{\infty}(M,E))$) the space of all bounded
operators from $\cD'(M,E)$ to $\cH^{\infty}$ (resp. from
$\cH^{-\infty}$ to $C^{\infty}(M,E)$). Since any operator from
$\Psi^{-N}(M,E)$ with $N>\dim M$ is a trace class operator in
$L^2(M,E)$, one can easily show the following inclusions
\begin{equation}\label{e:inclusions}
\begin{aligned}
\cL(\cD'(M,E),\cH^{\infty})&\subset
\cL^1(\cH^{-\infty},\cH^{\infty}),\\
\cL(\cH^{-\infty},C^{\infty}(M,E))&\subset
\cL^1(\cH^{-\infty},\cH^{\infty}).
\end{aligned}
\end{equation}

\begin{thm}\label{t:powers}
For any $\alpha\in \RR$, the operator $\langle
D\rangle^{\alpha}=(D^2+I)^{\alpha/2}$ can be written as $$\langle
D\rangle^{\alpha}=P(\alpha)+R(\alpha),$$ where:
\medskip\par
$($a$)$ $P(\alpha)\in \Psi^\alpha(M,E)$ is a self-adjoint, elliptic
operator with the positive, scalar principal symbol and the
holonomy invariant transversal principal symbol;

$($b$)$ For any $K\in \Psi^{*,-\infty}(M,{\mathcal F},E)$,
$KR(\alpha)\in\cL(\cH^{-\infty},C^{\infty}(M,E))$, and
$R(\alpha)K\in \cL(\cD'(M,E),\cH^{\infty})$.
\end{thm}
\begin{proof}
Using the standard construction of parametrix for elliptic
operators in some conic neighborhood of $N^*\cF$, one gets an
analytic family $C_1(\lambda),\lambda\not\in\RR_+$, of operators
from $\Psi^{-2}(M,E)$ such that
\begin{equation}\label{e:param}
C_1(\lambda)(D^2+I-\lambda I)=I-r_1(\lambda),\quad
\lambda\not\in\RR_+,
\end{equation}
where $r_1(\lambda)\in \Psi^{0}(M,E)$ has order $-\infty$  in some
conic neighborhood of $N^*\cF$ (see \cite{noncom} for more
details). Hence, we have $$ (D^2+I-\lambda
I)^{-1}=C_1(\lambda)+r_1(\lambda)(D^2+I-\lambda I)^{-1}, \quad
\lambda\not\in\RR_+. $$ Using the Cauchy integral formula with an
appropriate contour $\Gamma$ in the complex plane, we get  $$
(D^2+I)^{\alpha/2}=\frac{i}{2\pi}\int_{\Gamma}
\lambda^{\alpha/2-N}\langle D\rangle^{2N} (D^2+I-\lambda
I)^{-1}d\lambda =P_1(\alpha)+R_1(\alpha), $$ with some natural $N$
such that ${\rm Re}\;\alpha<2N$, where
\begin{gather*}
P_1(\alpha)=\frac{i}{2\pi}\int_{\Gamma}
\lambda^{\alpha/2-N}\langle D\rangle^{2N} C_1(\lambda)d\lambda, \\
R_1(\alpha)=\frac{i}{2\pi}\int_{\Gamma}
\lambda^{\alpha/2-N}\langle D\rangle^{2N}
r_1(\lambda)(D^2+I-\lambda I)^{-1}d\lambda.
\end{gather*}

In a standard manner (see \cite{noncom}), one can prove that
$P_1(\alpha)$ is a transversally elliptic operator of class
$\Psi^{\alpha}(M,E)$ with the scalar principal symbol and the
holonomy invariant, positive transversal principal symbol.

Let $K\in \Psi^{*,-\infty}(M,{\mathcal F},E)$. For any real $s$,
one can write $$ KR_1(\alpha)\langle
D\rangle^s=\frac{i}{2\pi}\int_{\Gamma} \lambda^{\alpha/2-N}
K\langle D\rangle^{2N}r_1(\lambda)\langle D\rangle^s(D^2+I-\lambda
I)^{-1}d\lambda. $$ By Lemma~\ref{smooth}, the operator $K\langle
D\rangle^{2N}r_1(\lambda)$ has the smooth kernel and defines a
bounded operator from $\cH^{-\infty}\subset \cD'(M,E)$ to
$C^{\infty}(M,E)$. Since $\langle D\rangle^s(D^2+I-\lambda
I)^{-1}$ maps $\cH^{-\infty}$ to $\cH^{-\infty}$, this implies
that the operator $KR_1(\alpha)$ is an operator of class
$\cL(\cH^{-\infty},C^{\infty}(M,E))$.

Taking adjoints in \eqref{e:param}, we get $$ (D^2+I-\lambda
I)C^*_1(\lambda)=I-r^*_1(\lambda),\quad \lambda\not\in\RR_+. $$ It
follows that
$C_1(\lambda)-C^*_1(\lambda)=C_1(\lambda)r^*_1(\lambda)-
r_1(\lambda)C^*_1(\lambda)$ has order $-\infty$ in some conic
neighborhood of $N^{*}\cF$. Moreover, using the formula
\begin{multline*}
\langle D\rangle^{2N}C_1(\lambda)-C_1^*(\lambda)\langle
D\rangle^{2N}\\
\begin{aligned}
=&\frac{1}{\lambda}\langle D\rangle^{2}\left(\langle
D\rangle^{2(N-1)}C_1(\lambda)- C_1^*(\lambda)\langle
D\rangle^{2(N-1)}\right)\langle D\rangle^{2}\\
&+\frac{1}{\lambda}\left(\langle D\rangle^{2(N-1)}r_1(\lambda)-
r_1^*(\lambda)\langle D\rangle^{2(N-1)}\right),
\end{aligned}
\end{multline*}
one can prove by induction that $\langle
D\rangle^{2N}C_1(\lambda)-C_1^*(\lambda)\langle D\rangle^{2N}$ has
order $-\infty$ in some conic neighborhood of $N^{*}\cF$. This
implies that the same is true for
$P_1(\alpha)-P^*_1(\alpha)=R^*_1(\alpha)-R_1(\alpha)$. Combining
Lemma~\ref{smooth} and duality arguments, we get that, for any
$K\in \Psi^{*,-\infty}(M,{\mathcal F},E)$, the operator
$R_1(\alpha)K=(K^*R_1(\alpha)+K^*(R^*_1(\alpha)-R_1(\alpha)))^*$
extends to a bounded operator from $\cD'(M,E)$ to $\cH^{\infty}$.

Let $P(\alpha)\in \Psi^{\alpha}(M,E)$ be a self-adjoint, elliptic
operator with the positive scalar principal symbol such that the
operator $P_1(\alpha)-P(\alpha)$ has order $-\infty $ in some
neighborhood of $N^*{\mathcal F}$ (see also \cite{dgtrace}) and
$R(\alpha)=\langle D\rangle^{\alpha/2}-P(\alpha)$. By
Lemma~\ref{smooth}, for any $K\in \Psi^{*,-\infty}(M,{\mathcal
F},E)$, the operator $K(P(\alpha)-P_1(\alpha))$ is a smoothing
operator, that immediately completes the proof.
\end{proof}

Let $\langle D\rangle=P+R$ be a representation given by
Theorem~\ref{t:powers}. Denote by $e^{itP}$ the strongly
continuous group of bounded operators in $L^2(M,E)$ generated by
the elliptic operator $i P$. Put also $R(t)=e^{it\langle
D\rangle}-e^{itP}$.

\begin{proposition}
\label{mapping} For any $K\in \Psi^{*,-\infty}(M,{\mathcal F},E)$,
$KR(t), t\in \RR,$ is a smooth family of operators from
$\cL(\cH^{-\infty},C^{\infty}(M,E))$, and $R(t)K, t\in \RR,$ is a
smooth family of operators from $\cL(\cD'(M,E),\cH^{\infty})$.
\end{proposition}

\begin{proof} By the Duhamel formula, for any $K\in \Psi^{*,-\infty}(M,{\mathcal
F},E)$ and $u\in H^1(M,E)\subset \Dom (P)$, one can write
\begin{equation*}
KR(t)u=i \int_0^t e^{i\tau P}e^{-i\tau P}Ke^{i\tau P}
         \,R\,e^{i(t-\tau)\langle D\rangle}u\,d\tau.
\end{equation*}
By Theorem~\ref{Egorov:elliptic}, $e^{-i\tau P}Ke^{i\tau P}\in
\Psi^{*,-\infty}(M,{\mathcal F},E)$. Therefore, the operator
$e^{-i\tau P}Ke^{i\tau P}R$ extends to a bounded operator from
$\cH^{-\infty}$ to $C^{\infty}(M,E)$. Since $e^{i\tau P}$ maps
$C^{\infty}(M,E)$ to $C^{\infty}(M,E)$ and $e^{i(t-\tau)\langle
D\rangle}$ is a bounded operator in $\cH^{-\infty}$, the operator
$KR(t)$ extends to a bounded operator from $\cH^{-\infty}$ to
$C^{\infty}(M,E)$.

Using the formula
\begin{equation}\label{e:induction0}
\frac{d^n}{dt^n}KR(t)=i\frac{d^{n-1}}{dt^{n-1}}KPR(t)+i^{n}KR\langle
D\rangle^{n-1}e^{it\langle D\rangle}, \quad n\in\NN,
\end{equation}
one can show by induction that, for any $K\in
\Psi^{*,-\infty}(M,{\mathcal F},E)$, the function $KR(t)$ is
smooth as a function on $\RR$ with values in
$\cL(\cH^{-\infty},C^{\infty}(M,E))$. The similar statement,
concerning to the operator $R(t)K$, follows by duality.
\end{proof}

\begin{proof}[Proof of Theorem~\ref{Egorov}] Let $\langle
D\rangle=P+R$ be a representation given by Theorem~\ref{t:powers}.
Let $K\in \Psi^{m,-\infty}(M,{\mathcal F},E)$. By
Theorem~\ref{Egorov:elliptic}, it follows that the operator
$\Phi_t^{P}(K)=e^{itP}Ke^{-itP}$ is in
$\Psi^{m,-\infty}(M,{\mathcal F},E)$. Moreover, if $E$ is the
trivial line bundle, the subprincipal symbol of $D^2$ vanishes,
and $k\in S^{m}(G_{{\mathcal F}_N},|T{\mathcal G}_N|^{1/2})$ is
the principal symbol of $K$, then the principal symbol $k(t)\in
S^m(G_{{\mathcal F}_N}, |T{\mathcal G}_N|^{1/2})$ of
$\Phi_t^{P}(K)$ is given by $k(t)=F^*_t(k)$.

To complete the proof, it suffices to show that
$\Phi_t(K)-\Phi_t^{P}(K), t\in \RR,$ is a smooth family of
operators of class $\cL^1(\cH^{-\infty},\cH^{\infty})$. We have $$
\Phi_t(K)-\Phi_t^{P}(K)=e^{itP}KR(-t)+R(t)Ke^{-it\langle
D\rangle}. $$ Using Proposition~\ref{mapping}, the fact that the
operator $e^{itP}$ takes $C^{\infty}(M,E)$ to itself and
(\ref{e:inclusions}), we get that $e^{itP}KR(-t)$ belongs to
$\cL^1(\cH^{-\infty},\cH^{\infty})$. Similarly, using
Proposition~\ref{mapping}, the fact that the operator
$e^{-it\langle D\rangle}$ is bounded in $\cH^{-\infty}$, and
(\ref{e:inclusions}), we get that
$R(t)K\in\cL(\cD'(M,E),\cH^{\infty}) \subset
\cL^1(\cH^{-\infty},\cH^{\infty})$ and, furthermore,
$R(t)Ke^{-it\langle D\rangle}\in
\cL^1(\cH^{-\infty},\cH^{\infty})$.
\end{proof}

\section{Noncommutative geometry of foliations}
Let $({\mathcal A},{\mathcal H},D)$ be a spectral triple
associated with a compact foliated Riemannian manifold
$(M,{\mathcal F})$ as in Section~\ref{ncg}. In this section, we
give a description of all the objects introduced in
Section~\ref{ncg} for this spectral triple. In particular, we will
prove Theorem~\ref{noncom:flow}.

First, we introduce a notion of scalar principal symbol for an
operator of class $\Psi^{m,-\infty}(M,{\mathcal F},E)$. Recall
that the bundle $\pi^*E$ on $N^*\cF$ is holonomy equivariant.
Therefore, there is a canonical embedding $$
i:C^{\infty}_{prop}(G_{{\mathcal F}_N},|T{\mathcal
G}_N|^{1/2})\hookrightarrow C^{\infty}_{prop}(G_{{\mathcal
F}_N},{\mathcal L}(\pi^*E)\otimes |T{\mathcal G}_N|^{1/2}),$$
which takes $k\in C^{\infty}_{prop}(G_{{\mathcal F}_N},|T{\mathcal
G}_N|^{1/2})$ to $i(k)=k\,\pi^*T$. We will identify
$C^{\infty}_{prop}(G_{{\mathcal F}_N},|T{\mathcal G}_N|^{1/2})$
with its image
\[
i(C^{\infty}_{prop}(G_{{\mathcal F}_N},|T{\mathcal
G}_N|^{1/2}))\subset C^{\infty}_{prop}(G_{{\mathcal
F}_N},{\mathcal L}(\pi^*E)\otimes |T{\mathcal G}_N|^{1/2}).
\]

We say that $P\in \Psi^{m,-\infty}(M,{\mathcal F},E)$ has the
scalar principal symbol if its principal symbol belongs to
$C^{\infty}_{prop}(G_{{\mathcal F}_N},|T{\mathcal G}_N|^{1/2})$.
Let $\Psi_{sc}^{m,-\infty}(M,{\mathcal F},E)$ denote the set of
all $K\in \Psi^{m,-\infty}(M,{\mathcal F},E)$ with the scalar
principal symbol. For any $k\in C^{\infty}_c(G,|T{\cG}|^{1/2})$,
the operator $R_E(k)$ is in $\Psi^{0,-\infty}_{sc}(M,\cF,E)$ and
its principal symbol $\sigma(R_E(k))$ is equal to $\pi_G^*k\in
C^{\infty}_{prop}(G_{{\mathcal F}_N},|T{\mathcal G}_N|^{1/2})$
where $\pi_G :G_{\cF_N}\to G$ is defined in Section~\ref{s:trpdo}.

Recall that $\delta$ denotes the inner derivation on ${\mathcal
L}({\mathcal H})$ defined by $\langle D\rangle$ (see
(\ref{e:delta})). It is easy to see that the class
$\cL^1(\cH^{-\infty},\cH^{\infty})$ belongs to the domain of
$\delta $ and is invariant under the action of $\delta$. Moreover,
one can easily show that $\cL^1(\cH^{-\infty}, \cH^{\infty})$ is
an ideal in ${\rm OP}^{0}$.

\begin{prop}\label{p:psi}
Any operator $K\in\Psi^{0,-\infty}(M,\cF,E)$ belongs to ${\rm
OP}_0^{0}$. Moreover, for any natural $n$ and for any
$K\in\Psi^{0,-\infty}(M,\cF,E)$, the operator $\delta^n(K)$
belongs to
$\Psi^{0,-\infty}(M,\cF,E)+\cL^1(\cH^{-\infty},\cH^{\infty})$. If
$K\in\Psi_{sc}^{0,-\infty}(M,\cF,E)$, $\delta^n(K)$ belongs to
$\Psi_{sc}^{0,-\infty}(M,\cF,E)+\cL^1(\cH^{-\infty},\cH^{\infty})$.
\end{prop}

\begin{proof} Let $\langle D\rangle=P+R$ be a representation given
by Theorem~\ref{t:powers}. Let $\delta_0$ denote the inner
derivation on ${\mathcal L}({\mathcal H})$ defined by $P$: $$
\delta_0(T)=[P,T], \quad T\in \Dom\delta_0\subset {\mathcal
L}({\mathcal H}). $$ Let $K\in \Psi^{0,-\infty}(M,\cF,E)$. Since
the principal symbol of $P$ is scalar and its transversal
principal symbol  is holonomy invariant, it is easy to see that
$\delta_0(K)$ is an operator of class $\Psi^{0,-\infty}(M,\cF,E)$,
that implies that $K$ belongs to the domain of $\delta^n_0$ for
any natural $n$.

We will prove by induction on $n$ that any $K\in
\Psi^{0,-\infty}(M,\cF,E)$ belongs to the domain of $\delta^n$ for
any natural $n$, and $$ \delta^n(K)-\delta^n_0(K)\in
\cL^1(\cH^{-\infty},\cH^{\infty}). $$ By Theorem~\ref{t:powers}
and (\ref{e:inclusions}), it follows that
\[
\delta(K)-\delta_0(K)=RK-KR\in \cL^1(\cH^{-\infty},\cH^{\infty}).
\]
Now assume that the statement holds for some natural $n$. Then one
can write
\begin{equation*}
\delta^{n+1}(K)-\delta^{n+1}_0(K) =
\delta(\delta^n(K)-\delta^n_0(K))+R\delta^n_0(K)-\delta^n_0(K)R,
\end{equation*}
that belongs to $\cL^1(\cH^{-\infty},\cH^{\infty})$, since
$\delta$ takes $\cL^1(\cH^{-\infty},\cH^{\infty})$ to itself and,
by Theorem~\ref{t:powers}, $R\delta^n_0(K)$ and $\delta^n_0(K)R$
are in $\cL^1(\cH^{-\infty},\cH^{\infty})$.

It remains to note that, by \cite{noncom}, for any
$K\in\Psi^{0,-\infty}(M,\cF,E)$, the operators $K\langle
D\rangle^{-1}$ and $\langle D\rangle^{-1}K$ are compact operators
in $L^2(M,E)$.
\end{proof}

Since $\cA = \cinf_c(G,|T\cG|^{1/2}) \subset
\Psi^{0,-\infty}(M,\cF,E)$, Proposition~\ref{p:psi} easily implies
the following

\begin{thm}\label{p:smooth}
For any $a\in {\mathcal A}$, the operators $a$ and $[D,a]$ belong
to ${\rm OP}^{0}$. Moreover, the algebra ${\mathcal B}$ generated
by $\delta^n(a), a\in {\mathcal A}, n\in {\NN}$ is contained in
${\rm OP}_0^{0}$.
\end{thm}

By Theorem~\ref{p:smooth}, it follows that the spectral triple
$({\mathcal A},{\mathcal H},D)$ is smooth. Next we will give a
description of $\cB$ and $\Psi^*_0({\mathcal A})$.

\begin{prop}\label{p:calB}
Any element $b\in \cB$ can written as $b=B+T$, where $B\in
\Psi_{sc}^{0,-\infty}(M,\cF,E)$ and $T\in
\cL^1(\cH^{-\infty},\cH^{\infty})$.
\end{prop}
\begin{proof}
By Proposition~\ref{p:psi}, the statement holds for any $b$ of the
form $\delta^n(a), a\in {\mathcal A}, n\in \NN$. Since
$\cL^1(\cH^{-\infty},\cH^{\infty})$ is an ideal in ${\rm OP}^{0}$,
this implies the statement for an arbitrary element of $\cB$.
\end{proof}

\begin{prop}\label{p:pseudo}
For any natural $N$, the algebra $\Psi^*_0({\mathcal A})$ is
contained in $\Psi_{sc}^{*,-\infty}(M,{\mathcal F},E)+{\rm
OP}_0^{-N}$.
\end{prop}

\begin{proof}
Take any $P\in \Psi^*_0({\mathcal A})$ of the form $P\sim
\sum_{j=0}^{+\infty}b_{q-j}\langle D\rangle^{q-j}$ with
$b_{q-j}\in {\mathcal B}$. Fix an arbitrary integer $j$. Let
$\langle D\rangle^j=P(j)+R(j)$ be a representation given by
Theorem~\ref{t:powers}. By Proposition~\ref{p:calB}, one can write
$b_j=B_j+T_j$, where $B_j\in \Psi_{sc}^{0,-\infty}(M,\cF,E)$ and
$T_j\in \cL^1(\cH^{-\infty},\cH^{\infty})$. So we have $$
b_j\langle D\rangle^j=B_jP(j)+B_jR(j)+T_j\langle D\rangle^j. $$
Here $B_jP(j)\in \Psi_{sc}^{j,-\infty}(M,{\mathcal F},E)$ (see
\cite{noncom}), $B_jR(j)\in \cL^1(\cH^{-\infty}, \cH^{\infty})$ by
Theorem~\ref{t:powers} and $T_j\langle D\rangle^j\in
\cL^1(\cH^{-\infty},\cH^{\infty})$ by the definition of
$\cL^1(\cH^{-\infty},\cH^{\infty})$. Thus, $b_j\langle
D\rangle^j\in \Psi^{j,-\infty}_{sc}(M,{\mathcal F},E) +
\cL^1(\cH^{-\infty}, \cH^{\infty})$, that completes the proof.
\end{proof}

Now we need the following result on continuity of the principal
symbol map given by (\ref{e:symbol}). Let $E$ be a vector bundle
on a compact foliated manifold $(M,\cF)$. Denote by
$\bar{\Psi}^{0,-\infty}(M,{\mathcal F},E)$ the closure of
$\Psi^{0,-\infty}(M,{\mathcal F},E)$ in the uniform topology of
${\mathcal L}(L^2(M,E))$.

\begin{prop}
\label{seq} $($1$)$ The principal symbol map $$ \sigma :
\Psi^{0,-\infty}(M,{\mathcal F},E)\rightarrow S^{0}(G_{{\mathcal
F}_N},{\mathcal L}(\pi^*E)\otimes |T{\mathcal G}_N|^{1/2}) $$
extends by continuity to a homomorphism $$
\bar{\sigma}:\bar{\Psi}^{0,-\infty}(M,{\mathcal F},E) \rightarrow
\bar{S}^{0}(G_{{\mathcal F}_N},{\mathcal L}(\pi^*E)\otimes
|T{\mathcal G}_N|^{1/2}). $$

$($2$)$ The ideal $I_\sigma=\Ker \bar{\sigma}$ contains the ideal
${\mathcal K}$ of compact operators in $L^2(M,E)$.
 \end{prop}

Proposition~\ref{seq} can be proven by an easy adaptation of the
proof of analogous fact for pseudodifferential operators on
compact manifolds (see, for instance, \cite{Palais,Seeley65}).
\bigskip
\begin{proof}[Proof of Theorem~\ref{noncom:flow}] By
Proposition~\ref{p:pseudo}, it follows that the algebra ${\mathcal
C}_0$ is contained in $\Psi^{0,-\infty}(M,{\mathcal F})+{\rm
OP}^{-N} (\cH^{-\infty},\cH^{\infty})$ for any $N$ and its
closure, $\bar{\mathcal C}_0$, is contained in
$\bar{\Psi}^{0,-\infty}(M,{\mathcal F})+\cK$. By
Proposition~\ref{seq}, the principal symbol map $\bar{\sigma}$
induces a map $P: S^*{\mathcal A}\rightarrow
\bar{S}^{0}(G_{{\mathcal F}_N},|T{\mathcal G}_N|^{1/2}).$ By
Theorem~\ref{Egorov}, it follows that the diagram \eqref{e:cd} is
commutative that completes the proof.
\end{proof}

\begin{rem}
\label{extens} Suppose $E$ is a holonomy equivariant vector
bundle. Let $C^{*}_E(G)$ be the closure of
$R_E(\cinf_c(G,|T\cG|^{1/2}))$ in the uniform operator topology of
${\mathcal L}(L^2(M,E))$ and $C^{\ast}_{r}(G)$ the reduced
foliation $C^*$-algebra (see, for instance, \cite{F-Sk}). By
\cite{F-Sk}, there is a natural surjective projection $ \pi_E :
C^{\ast}_{E}(G) \rightarrow C^{\ast}_{r}(G).$ The map $\pi_G
:G_{\cF_N}\to G$ defines a natural embedding
$C^{\ast}_{r}(G)\subset \bar{S}^{0}(G_{{\mathcal F}_N},
|T{\mathcal G}_N|^{1/2}).$ Since $R_E(k)\in
\Psi^{0,-\infty}(M,{\mathcal F},E)$ for any
$\cinf_c(G,|T\cG|^{1/2})$, $C^{\ast}_{E}(G)$ is contained in
$\bar{\Psi}^{0,-\infty}(M,{\mathcal F},E)$. Moreover, the
restriction of $\bar{\sigma}$ to $C^{\ast}_{E}(G)$ coincides with
$\pi_E$. So the principal symbol map $\bar{\sigma}$ provides an
extension of $\pi_E$ to $\bar{\Psi}^{0,-\infty}(M,{\mathcal
F},E)$. In particular, if $I_\sigma=\Ker \bar{\sigma}$ coincides
with ${\mathcal K}$, then $\pi_E$ is injective, and the holonomy
groupoid $G$ is amenable (see, for instance, \cite{claire-jean}).
\end{rem}
\hyphenation{grou-po-ids}

{\bf Acknowledgements.} The author acknowledges hospitality and
support of the Ohio State University where the work was completed
as well as partial support from the Russian Foundation for Basic
Research, grant no. 04-01-00190. We also thank the referees for
corrections and suggestions.

\end{document}